\documentclass[a4paper,10pt]{article}

\newtheorem{thm}{Theorem}
\newtheorem{cor}[thm]{Corollary}
\newtheorem{rem}[thm]{Remark}
\newtheorem{lem}[thm]{Lemma}
\newtheorem{prop}[thm]{Proposition}
\def \a{{\alpha}}
\def \b{{\beta}}
\def \D{{\Delta}}
\def \DD{{\mathcal D}}

\def \e{{\varepsilon}}

\def \g{{\gamma}}
\def \G{{\Gamma}}

\def \O{{\Omega}}

\def \t{{\vartheta}}
\def \ttau{{\theta}}
\def \m{{\mu}}
\def \s{{\sigma}}

\def \A{{\cal A}}

\def \E{{\bf E}\, }

\def \NNN{{\cal N}}

\def \P{{\bf P}}
\def \q{{\quad}}

\def \qq{{\qquad}}
\def \R{{\bf R}}

\def \T{{\bf T}}
\def \td{{\widetilde D}}
\def \ua{{\underline{a}}}

\def \ud{{\underline{d}}}

\def \uz{{\underline{z}}}

\def \Z{{\bf Z}}
\def \zz{{\cal Z}}

\def \noi{{\noindent}}

\def\qed{\hbox{\vrule height 6pt depth 0pt width
6pt}}
\def\cqfd{\hfill\penalty 500\kern 10pt\qed\medbreak}
\font\phh=cmcsc10
at  8 pt

\title{On the Supremum of
Some
Random Dirichlet Polynomials
}
\author{Mikhail Lifshits and Michel Weber}

\begin{document}

\maketitle

\begin{abstract}
We study the average supremum  of some random Dirichlet polynomials
$D_N(t)=\sum_{n=1}^N\e_n d(n) n^{-\s - it}$,
where $(\e_n)$ is a sequence of independent Rademacher random variables,
the weights $(d(n))$ satisfy some reasonable conditions and $0\le \s \le 1/2$.
We use an approach
  based on  methods of stochastic processes, in particular the metric
entropy method developed in \cite{LW}.
\end{abstract}


\section{Introduction}
Let $\{d(n), n\ge 1\}$ be a sequence of real numbers.   Consider the Dirichlet polynomials
$$ P( \s + it)=P_N( \s + it)= \sum_{n=1}^N  d(n) n^{ -\s - it}  .
$$As is   well-known,   the abscissa of uniform convergence $\s_u$ of the associated Dirichlet
series $\sum_{n=1}^\infty  d(n) n^{ -\s - it}$, which is  defined by
$$
\s_u=\inf\Big\{\s : \sum_{n=1}^\infty  d(n) n^{-\s - it}\
\hbox{converges uniformly over $t\in\R$} \Big\},
$$
satisfies the relation
 $\s_u= \limsup_{N\to \infty} {{\log  \,  \sup_{t\in \R} |P_N(it)  |} \over \log N}   .
 $
And this    gives  a motivation  to study of the supremum of the Dirichlet polynomials
$ P_N(s)  $
over lines $\{s=\s+it, \ t\in \R\}$  (see for instance   \cite{KQ}, \cite{LW} and the references therein).
 \par
\medskip\par A first basic reduction step  allows to replace the
Dirichlet polynomial by some relevant   trigonometric
polynomial. Introduce some
necessary notation. Let $2=p_1<p_2<\ldots$ be the sequence of all
primes. Let $\pi(N)$ denote  the number of prime numbers
that are less or equal to $N$. Now fix $N$ and  put $\tau =\pi(
N)$. If $  n=\prod_{j=1}^\tau p_j^{a_j(n)}$, we write
$\ua(n)=\big\{a_j(n), 1\le j\le \tau\big\}$.   Let also $\T=[0,1[=\R/\Z$ be the torus.
Define  for $\uz= (z_1,\ldots, z_\tau) \in \T^\tau$
$$
 Q(\uz)= \sum_{n=1}^N d(n) n^{-\s}e^{2i\pi\langle \ua(n),\uz\rangle}.
$$
  H. Bohr's observation (see e.g. \cite{Q1}), based on Kronecker's Theorem (see \cite{HW}, Theorem 442, p.382) states that
\begin{equation} \label{e11}
 \sup_{t\in \R} \big|P(\s +it)\big| =\sup_{\uz \in \T^\tau} \big|Q(\uz)\big|\ .
\end{equation}

The supremum properties of random  Dirichlet polynomials and   random
Dirichlet series were investigated by   Hal\'asz, Bayard, Konyagin, Queff\'elec,   recently  by the authors,
 and earlier, by Hartman, Clarke, Dvoretzky and Erd\" os, where random
power  series are also considered  (see \cite{LW} for references).
 
Let $\e= \{\e_n, n\ge 1\}$  be (here and throughout the whole
paper) a sequence of independent Rademacher random variables
($\P\{ \e_i=\pm 1\} =1/2$) defined on a basic
probability space $(\O, \A, \P)$.
Consider the random Dirichlet polynomials
\begin{equation} \label{e12}
  \DD(s) =\sum_{n=1}^N \e_n d(n) n^{-s}.
\end{equation}

In the particular case $d(n)\equiv 1$,
if $\s=0$,    the following result was proved by Hal\'asz (see \cite{Q2},\cite{Q3}): for some absolute
constant $C$, and all integers $N\ge 2$

\begin{equation} \label{e13}
C ^{-1} {N  \over   \log  N}
\le
\E\, \sup_{t \in \R} \left|\sum_{n=1}^N \e_n  n^{  - it} \right|
\le
C {N\over \log N}\ .
\end{equation}

 In \cite{Q2},\cite{Q3}
(see also \cite{Q1} for a first result), Queff\'elec extended Hal\'asz's
result to the range of values $0\le \s<1/2$; and provided a probabilistic
proof of the original one, using Bernstein's inequality for polynomials,
properties of complex Gaussian processes and the sieve method introduced by
Hal\'asz. He obtained that for some constant $C_\s$ depending on $\s$
only, and all integers $N\ge 2$
\begin{equation} \label{e14}
C_\s ^{-1} {N^{1-\s}  \over   \log  N }
\le
\E\, \sup_{t \in \R}\left|\sum_{n=1}^N \e_n  n^{-\s  - it}\right|
\le C_\s {N^{1-\s} \over \log  N }\ .
\end{equation}
This result was further extended in \cite{LW} for the polynomials
\[
\sum_{{P^+(n) \le \ttau,\atop n\le N}} \e_n  n^{-\s-it}
\]
(here and in what follows  $P^+(n)$ denotes the largest prime divisor
of $n$) with fine estimates depending
on both parameters $\ttau$ and $N$. For small $\ttau$ the estimates of such kind are related
to the construction of so called Rudin-Shapiro Dirichlet polynomials.
 \bigskip\par
In this work, we show that the approach developed in \cite{LW} is sufficiently
robust to allow the similar study for random Dirichlet polynomials with
reasonable weights. We will also extend in a separate section the main result in
\cite{LW} to the boundary case $\s=1/2$.

We introduce now some characteristics of weights. Let
\[
 D_1(M)= \sum_{m=1}^M d(m); \qquad  \td_1(M)=\max_{1\le m\le M} \frac {D_1(m)}{m} \ ,
\]
and, similarly,
\[
 D_2(M)= \sum_{m=1}^M d(m)^2; \qquad  \td_2^2(M)=\max_{1\le m\le M} \frac {D_2(m)}{m} \ ,
\]
Obviously, we have $\td_2(M)\ge \td_1(M)$.
\medskip

{\bf Example 1.}\ Let $d(n)$ be the divisor function (number of divisors of integer $n$).
Then $d$ is a multiplicative function and it is well-known that
\[ \td_1(M) \sim \log M, \qquad  \td_2(M) \sim \log^{3/2} M.
\]
\medskip

{\bf Example 2.}\  Let von Mangoldt function $\Lambda$ be defined as follows:

\begin{equation} \label{mangoldt}
\Lambda(n)=
\cases{\log p,&\ $n=p^k$,\ $p$\ \textrm{is prime} \cr
       0,     & \textrm{else} .
      }
\end{equation}

Then  $\Lambda$ is neither additive nor multiplicative, and we have for any integers $k,j$
and any prime $p$ either $\Lambda(k p^j)=0$, or $k=p^i$, hence
$\Lambda(k p^j)=\Lambda(k)=\log p$.

Moreover, elementary calculations show that
\[ \lim_{M\to\infty} \td_1(M) = 1, \qquad  \td_2(M) \sim \log M.
\]
\medskip

\begin{thm} \label{t1}
 Let $0\le\s \le 1/2$ and assume that
\begin{equation} \label{basic_d}
      d(k p^j) \le C d(k) j^H
\end{equation}
 for some positive $C,H$, any positive integer $k,j$ and any prime $p$.
 Then there exists a constant $C$ depending on $d$ and $\s$ such that for any integer $N\ge 2$
it is true that
\begin{equation}\label{main_upper}
\E \sup_{t\in \R} \big|\DD(\s+it) \big| \le {C\ N^{1-\s}  \td_2(N)  \over \log N} \ .
\end{equation}
Moreover, if for some $b<b_*:=(\sqrt{5}-1)/4\approx 0.31$
\begin{equation}\label{extra}
\td_2(M) \le C  M^b,
\end{equation}
then
\[
\E \sup_{t\in \R} \big| \DD(\s+it) \big| \le {C\ N^{1-\s}   \over \log N} \ .
\]
\end{thm}

\begin{rem} \nonumber \rm
  If one assumes $(\ref{extra})$ with $b\in(b_*,1/2)$, then our proof shows that
\[
\E \sup_{t\in \R} \big| \DD(\s+it) \big| \le {C\ N^{r-\s}
\over (\log N)^{1/2}}, \qquad \ r=r(b)=b+ \frac{2b+3}{4(1+b)},
\]
which is better than $(\ref{main_upper})$.
\end {rem}\smallskip \par

\begin{rem} \nonumber
\rm If $d(n)$ is a multiplicative function, then condition (6) is satisfied iff
$$ d(p^{r+j})\le C d(p^r) j^H,$$
for some $C>0$, $H>0$ and any $j\ge 1$, $r\ge 0$. This last condition is satisfied for
instance if $d(p)={\mathcal O}(1)$ and
$$\frac{d(p^{k+1})}{d(p^{k })} \le (1+\frac{1}{k})^H, \qquad  k=1,2, \ldots$$
A multiplicative function being completely defined by its values $ d(p^k)$, these ones
can be prescribed arbitrarily. This observation shows that condition (6) is satisfied for a
very large class of multiplicative functions.
\end {rem}


\section{Proof of Theorem \ref{t1}.}

For proving the upper bound, we first operate  the reduction to the study of a random
polynomial $Q$ on the multidimensional torus by using
(\ref{e11}). Next we use a decomposition $Q=Q_1+ Q_2$. The study of the  supremum
of the polynomial $Q_1$ is made by using the metric entropy method.
The investigation of the supremum of the polynomial $Q_2$ first relies upon
the contraction principle, reducing the study to the one of a
complex valued Gaussian process, next via
Slepian's comparison Lemma, by a thorough study of the $L^2$-metric
induced by this process.

  The calculations are first provided for $\s<1/2$.
At the end of the proof we comment on the case $\s=1/2$ for which some minor
details are different from the generic case.

Introduce some notation.
We can represent $[1,N]$ as the union of disjoint sets
$$
E_j=\big\{ 2\le n\le N :  P^+(n)=p_j\big\}, \q j=1,\ldots, \tau.
$$
For $\uz \in \T^\tau$ we put
$$
Q(\uz)= \sum_{j=1}^\tau \sum_{n\in E_j}
\e_n  d(n) n^{-\s}e^{2i\pi\langle \ua(n),\uz\rangle}.
$$
By  (\ref{e11}) we have
$$
\sup_{t\in \R}\big|\sum_{j=1}^\tau \sum_{n\in E_j} \e_n
d(n) n^{ -\s - it}\big| =\sup_{\uz \in \T^\tau}\big|Q(\uz)\big|.
$$
Let $1\le \nu <\tau$ be fixed. Write $Q=Q_1+Q_2$ where
\begin{eqnarray*}
Q_1(\uz) &=& \sum_{  P^+(n) \le p_{\nu }} \e_n d(n) n^{-\s}
e^{2i\pi\langle\ua(n),\uz\rangle}, \
\cr
Q_2(\uz) &=& \sum_{p_{\nu}<  P^+(n) \le  p_{\tau }}
\e_n  d(n) n^{-\s}e^{2i\pi\langle \ua(n),\uz\rangle}.
\end{eqnarray*}
First, evaluate the supremum of $Q_2$. Introduce the
following random process
$$
X^\e(\g) =\sum_{ \nu <j\le \tau} \a_j
\sum_{n\in E_j}  \e_n  d(n) n^{-\s} \b_{{n\over p_j }},
\qq \g\in \G,
$$
where
$\g =\big((\a_j)_{\nu <j\le \tau}, (\b_m)_{1\le m\le N/2}\big)
$
and
$$\G=\big\{ \g  :  |\a_j|\vee |\b_m|\le \!\!1,  \nu < \!j\le \tau,
1\! \le m\le N/2\big\}.
$$
Writing
\begin{eqnarray*}
Q_2(\uz)&=& \sum_{ \nu <j\le \tau}
e^{2i\pi z_j  }\sum_{n\in E_j}
\e_n d(n) n^{-\s} e^{2i\pi\{\sum_{k\not = j} a_k(n)z_k+ [a_j(n)-1]z_j\}}
\cr
&=& \sum_{ \nu <j\le \tau} e^{2i\pi z_j  }\sum_{n\in E_j}
\e_n  d(n)  n^{-\s}e^{2i\pi\left\{\sum_{k} a_k({n\over p_j}) z_k\right\}}
\end{eqnarray*}
and considering separately the imaginary and real parts of the exponents, it
follows that $Q_2(\uz)$ can be written as the sum of four terms each being of the form
$$
\eta\ \sum_{ \nu <j\le \tau}
\a_j \sum_{n\in E_j}  \e_n d(n) n^{-\s} \b_{{n\over p_j }},
$$
where $\eta\in\{1,i,-i,-1\}$, and
$$
\a_j= \a_j(\uz)=
\cases{  \cos (2\pi z_j), &\cr
 \hbox{\rm or}  &\cr
 \sin (2\pi z_j),&
}
\qquad \nu <j\le \tau;
$$
$$
\b_m = \b_m(\uz) =
\cases{ \cos\left( 2\pi \sum_{k} a_k({m}) z_k \right),  &\cr
\hbox{\rm or}  &\cr
\sin\left( 2\pi \sum_{k} a_k({m}) z_k \right),&
}
\qquad 1\le m\le {N\over 2}\ .
$$

Therefore, we obtain
$$\sup_{\uz \in \T^\tau}\big|Q_2(\uz)\big|
\le
4\sup_{\g \in\G}  \big|X^\e(\g)\big|.
$$
By the contraction principle (\cite{K} p.16-17)

$$
\E\, \sup_{\uz \in \T^\tau}\big|Q_2(\uz)\big|
\le  4 \ \sqrt{ \pi \over 2  }\
\E\, \sup_{\g \in \G}\big|X(\g)\big|,
$$
where $\{X(\g), \g\in \G\}$ is the same process as $X^\e(\g)$
except that the Rademacher random variables $\e_n$ are replaced by
independent
${\cal N}(0,1)$ random variables $\m_n$:
$$  X(\g) =\sum_{\nu<j\le \tau} \a_j \sum_{n\in E_j}
\m_n d(n) n^{-\s} \b_{{n\over p_j}}.
$$
The problem now reduces to estimating the
supremum of the real valued Gaussian process $X$. Towards  this aim,
we examine the $L^2$-norm of its increments:
\begin{eqnarray*}
\|X_\g-X_{\g'}\|_2^2
&=&
\sum_{ \nu <j\le \tau} \sum_{n\in E_j}
 d(n)^2 n^{-2\s} \big[\a_j\b_{{n\over p_j}}-\a'_j \b'_{{n\over p_j}}\big]^2
\cr
&\le&
2\!\sum_{ \nu <j\le \tau}\sum_{n\in E_j}  d(n)^2 n^{-2\s} \big[(\a_j -\a'_j)^2+
(\b_{{n\over p_j}}-\b'_{{n\over p_j}})^2\big]  ,
\end{eqnarray*}
where we have used the identity
\[
\a_j\b_{{n\over p_j}}-\a'_j \b'_{{n\over p_j}}= (\a_j -\a'_j)
\b_{{n\over p_j}}+ (\b_{{n\over p_j}}- \b'_{{n\over p_j}})\a'_j.
\]

\noi The "$\a$" component part is easily controlled as follows,
\begin{eqnarray}
\sum_{ \nu <j\le \tau} \sum_{n\in E_j}  d(n)^2 n^{-2\s} (\a_j -\a'_j)^2
&\le&
C\ \sum_{\nu<j\le\tau} (\a_j -\a'_j)^2p_j^{-2\s}
\sum_{m\le N/p_j}  \frac{d(m)^2}{ m^{2\s}}
\cr
\cr  \label{e21}
&\le&
C \sum_{\nu<j\le\tau} (\a_j -\a'_j)^2 {N^{1-2\s}\td _2^2(N/p_j) \over p_j},
\end{eqnarray}
since by Abel transformation
\begin{equation} \label{abel}
 \sum_{m\le M} d(m)^2 m^{-2\s} \le C \td_2^2(M) M^{1-2\s}.
\end{equation}

For the "$\b$" component part, we use $d(mp_j)\le C d(m)$,
which is a particular case of (\ref{basic_d}), and obtain
\begin{eqnarray}
\sum_{\nu <j\le \tau} \sum_{n\in E_j}
{   d(n)^2 (\b_{{n\over p_j}}- \b'_{{n\over p_j}})^2\over n^{2\s}}
&\le&
C \sum_{m \le N/p_\nu}
(\b_m-\b'_m)^2 \big(\sum_{{\nu<j\le \tau\atop  mp_j\le N}}
{ d(m)^2 \over (mp_j)^{2\s}} \big)
\cr                                     \label{e22}
&:=&
C \sum_{m \le N/p_\nu} K_m^2
(\b_m-\b'_m)^2.
\end{eqnarray}

Now we evaluate the coefficients $K_m$.

%
%
\smallskip\par
\noindent Take a unique $k\in (\nu,\tau]$ such that
$N/p_k< m \le N/p_{k-1}$. By using
\begin{equation} \label{e23}
  p_j \sim j \ \log j
\end{equation}
we have
\begin{eqnarray*}
K_m^2 &=& \sum_{\nu<j\le k-1 }  d(m)^{2} (mp_j)^{-2\s}
\le
d(m)^{2} m^{-2\s}\ \sum_{j\le k-1 } p_j^{-2\s}
\cr
&\le&
C\ d(m)^{2} m^{-2\s}\ \sum_{j\le k } (j\log j)^{-2\s}
\le
C\  d(m)^{2} m^{-2\s}\  {k^{1-2\s} \over (\log k)^{2\s}}
\cr
&\le&  C d(m)^{2} m^{-2\s} \  {k \over p_k^{2\s}}
 \le
C m^{-2\s}  d(m)^{2} \  {k \over (N/m)^{2\s}}
\cr
&=& C \  d(m)^{2} {k \over N^{2\s}}\ .
\end{eqnarray*}
Since $k\log k\le Cp_k\le C\ {N\over m}$, we have
$$ k\le C\ {N\over m} \ (\log({N\over m}))^{-1}.
$$
We arrive at
\begin{equation} \label{e23a}
K_m\le C\  d(m)   N^{-\s}   ({N\over m})^{1/2} \ (\log({N\over m}))^{-1/2}
\ .
\end{equation}
By Abel transformation,
\begin{eqnarray*}
& &\sum_{m\le M} ({N\over m})^{1/2} \ (\log({N\over m}))^{-1/2} d(m)
\cr
&\le& D_1(M) ({N\over M})^{1/2} \ (\log({N\over M}))^{-1/2}
+ C\ \sum_{m\le M} {N^{1/2}\over m^{3/2}} \ (\log({N\over m}))^{-1/2} D_1(m)
\cr
&\le& \td_1(M) \left( (N M)^{1/2} \ (\log({N\over M}))^{-1/2}
+ C\ \sum_{m\le M} {N^{1/2}\over m^{1/2}} \ (\log({N\over m}))^{-1/2} \right)
\cr
&\le& \td_1(M) \left( \frac{ (N M)^{1/2}}{ (\log({N\over M}))^{1/2}}
+ C N \ \int_0^{M/N} u^{-1/2} \ (\log(1/u))^{-1/2} du \right)
\cr
&\le& C \td_1(M)  (N M)^{1/2} \ (\log({N\over M}))^{-1/2}.
\end{eqnarray*}
It follows that
\begin{eqnarray}
\sum_{m\le N/p_\nu} K_m
&\le& C\  N^{-\s}\ \sum_{m\le N/p_\nu}
({N\over m})^{1/2} \ (\log({N\over m}))^{-1/2} d(m)
\cr \label{e23b}
&\le&
{C N^{1-\s} \td_1(N/p_\nu)\over \nu^{1/2} \log\nu}\ .
\end{eqnarray}

Now define a second Gaussian process by putting for all $\g\in \G$
$$
Y(\g) =
\sum_{\nu<j\le\tau} \big({\td_2^2(N/p_j) N^{1-2\s}\over p_j}\big)^{{1/ 2}}\a_j\xi'_j
+
\sum_{m\le N/p_\nu} K_m \ \b_m\xi''_m
 :=
 \ Y'_\g + Y''_\g ,
$$
where $\xi'_i  $, $\xi''_j$  are independent ${\cal N}(0,1)$ random
variables. It follows from (\ref{e21}) and  (\ref{e22}) that for some suitable constant
$C$, one has the comparison relations: for all $\g, \g'\in \G$,
$$
\|X_\g-X_{\g'}\|_2\le C \|Y_\g-Y_{\g'}\|_2.
$$
By virtue of the Slepian comparison lemma (see \cite{L}, Theorem 4 p.190), since
$X_0=Y_0=0$,
we have
$$
\E\, \sup_{\g\in \G}     |X_\g|\le
2 \E\, \sup_{\g\in \G}      X_\g \le
2 C \E\, \sup_{\g\in \G} Y_\g \le
2 C \E\, \sup_{\g\in \G} |Y_\g|.
$$
It remains to evaluate the supremum  of $Y$. First of all,
$$
\E\, \sup_{\g\in \G}|Y'(\g)|
\le N^{{1\over 2}- \s}\sum_{\nu<j\le\tau}   p_j^{-1/2} \td_2(N/p_j).
$$
By using (\ref{e23}) we have
$$
\sum_{\nu<j\le\tau}   p_j^{-1/2}
\le
   \sum_{1<j\le\tau}   p_j^{-1/2}
\le
 {C \tau^{1/2} \over (\log\tau)^{1/2}}\ ,
$$
thus
\begin{equation} \label{e24}
 \E\, \sup_{\g\in \G}|Y'(\g)|
 \le C\ N^{{1\over 2}-\s}  \td_2(N/p_\nu) \ {\tau^{1/2} \over (\log\tau)^{1/2}}
 \le {C\ N^{1-\s}  \td_2(N/p_\nu)  \over \log N}
 \ .
\end{equation}

Under assumption (\ref{extra}) we can get a better estimate by using
the following lemma.

\begin{lem} \label{lem2} Let $f(N)\le c N^b$, $b<1/2$. Then
\[
  \sum_{1\le j\le\tau(N)}   p_j^{-1/2} f(N/p_j)
\le
 {C N^{1/2} \over \log N}\ ,
\]
with $C$ depending on $c$ and $b$.
\end{lem}

{\bf Proof.}\
\begin{eqnarray*}
  \sum_{1\le j\le\tau(N)}   p_j^{-1/2} f(N/p_j) &\le&
 c\ \sum_{1\le j\le\tau(N)}   p_j^{-1/2} (N/p_j)^b \cr
 &=& c N^b\ \sum_{1\le j\le\tau(N)}
p_j^{-(1/2+b)}
 \le C N^b\ \sum_{1\le j\le\tau(N)}  \frac 1{(j\log j)^{1/2+b}}
 \cr
 &\le& C N^b\ \tau(N)^{1/2-b} (\log\tau(N))^{-1/2-b}
\le {C N^{1/2} \over \log N}.
\end{eqnarray*}
\cqfd

By applying lemma to $f(N)=\td_2(N)$ we see that (\ref{extra}) implies
\begin{equation} \label{e24_ex}
 \E\, \sup_{\g\in \G}|Y'(\g)|
 \le {C\ N^{1-\s}  \over \log N}
 \ .
\end{equation}

\noi To control the supremum of $Y''$, we use our
estimates for the sums of $K_m$ and write that
 \begin{equation} \label{e25}
 \E\, \sup_{\g\in \G} |Y''(\g)|
\le \sum_{m\le  N/p_\nu} K_m
\le
{C  N^{1-\s} \td_1(N/p_\nu)\over \nu^{1/2} \log\nu} \ .
\end{equation}
\bigskip

Now, we turn to the supremum of $Q_1 $. Introduce the auxiliary Gaussian process
$$
\Upsilon (\uz) =
\sum_{n: P^+(n) \le p_\nu } d(n)  n^{-\s}
\big\{\t_n \cos 2\pi \langle \ua(n),\uz\rangle +
\t_n'\sin 2\pi \langle \ua(n),\uz\rangle \big\}
,\qq  \uz\in \T^{\nu },
$$
where $\t_i$, $\t'_j$  are independent ${\cal N}(0,1)$
random variables. By symmetrization (see e.g. Lemma 2.3 p. 269 in \cite{PSW}),
\[
\displaystyle{\E\, \sup_{\uz \in \T^{\nu }}\big|Q_1(\uz)\big|\le
\sqrt{8\pi}
\E\, \sup_{\uz \in  \T^{\nu }}\big|\Upsilon (\uz)\big|},
\]
so that we are again led to evaluating the supremum of a real valued
Gaussian process. For  $\uz, \uz' \in \T^{\nu}$ put
$\big\|\Upsilon(\uz)-\Upsilon(\uz)\big\|_2 := \,\rho (\uz, \uz')$,
and observe that
\begin{eqnarray*}
& &\,\rho(\uz,\uz')^2
= 4 \sum_{n: P^+(n)\le p_\nu} {d(n)^2\over n^{2\s}} \sin^2(\pi \langle\ua(n),\uz -\uz'\rangle)
\cr
&\le&
4\pi^2\ \sum_{n: P^+(n) \le p_\nu}  {d(n)^2\over n^{2\s}} |\langle\ua(n),\uz -\uz'\rangle|^2
\cr
&\le& 4\pi^2\ \sum_{n: P^+(n)\le p_\nu}  d(n)^2 n^{-2\s}
   \left[ \sum_{j=1}^\nu a_j(n) |z_j - z'_j|\right]^2
\cr
&=& 4\pi^2\ \sum_{n: P^+(n)\le p_\nu}   \sum_{j_1,j_2=1}^\nu a_{j_1}(n)
    a_{j_2}(n) |z_{j_1} - z'_{j_1}|\ |z_{j_2} - z'_{j_2}|d(n)^2 n^{-2\s}
\cr
&=& 4\pi^2
\ \sum_{j_1,j_2=1}^\nu \sum_{n: P^+(n)\le p_\nu} a_{j_1}(n)
     a_{j_2}(n)  |z_{j_1} - z'_{j_1}|\ |z_{j_2} - z'_{j_2}| d(n)^2 n^{-2\s}
\cr
&\le& 4\pi^2\ \sum_{j_1,j_2=1}^\nu  |z_{j_1} - z'_{j_1}|\ |z_{j_2} -z'_{j_2}|
    \sum_{b_1,b_2=1}^\infty b_1 b_2 \sum_{{ n\le N, a_{j_1}(n)=b_1,\atop a_{j_2}(n)=b_2}} d(n)^2 n^{-2\s}
\cr
&\le&   C\sum_{j_1,j_2=1}^\nu |z_{j_1} - z'_{j_1}|\ |z_{j_2} - z'_{j_2}|
    \sum_{b_1,b_2=1}^\infty  \frac{b_1^{1+H} b_2^{1+H}} {p_{j_1}^{2 b_1\s} p_{j_2}^{2 b_2\s}}
\sum_{k\le N p_{j_1}^{-b_1} p_{j_2}^{-b_2}} d(k)^2 k^{-2\s}.
\end{eqnarray*}
By Abel transformation argument (\ref{abel}),
denoting $\D:=\td_2(N)$, we have
\begin{eqnarray*}
  \rho(\uz,\uz')^2
 &\le& C N^{1-2\s} \D^2
\sum_{j_1,j_2=1}^\nu  |z_{j_1} - z'_{j_1}|\ |z_{j_2} - z'_{j_2}|
\sum_{b_1,b_2=1}^\infty \frac{b_1^{1+H} b_2^{1+H}} {p_{j_1}^{2 b_1\s} p_{j_2}^{2 b_2\s}}
      [p_{j_1}^{-b_1} p_{j_2}^{-b_2}]^{1-2\s}
\cr
&=&  C N^{1-2\s} \D^2 \sum_{j_1,j_2=1}^\nu  |z_{j_1} - z'_{j_1}|\ |z_{j_2} - z'_{j_2}|
    \sum_{b_1,b_2=1}^\infty \frac{b_1^{1+H} b_2^{1+H}} {p_{j_1}^{b_1} p_{j_2}^{b_2}}
\cr
&=& C N^{1-2\s} \D^2  \left\{ \sum_{j=1}^\nu  |z_{j} - z'_{j}|
    \sum_{b=1}^\infty b^{1+H} \ p_{j}^{-b} \right\}^2.
\end{eqnarray*}
Thus,
 \begin{equation} \label{e27}
 \rho(z,z')\le C N^{1/2-\s} \D \left\{
   \sum_{j=1}^\nu  |z_{j} - z'_{j}|
   \sum_{b=1}^\infty  b^{1+H} \  p_{j}^{-b} \right\}.
\end{equation}
\medskip

Now we explore the entropy properties of the metric space $(\T^\nu,\rho)$. Towards
this aim, take $\e\in (0,1)$ and cover $T^\nu$ by rectangular cells so that if
$z$ and $z'$ belong to the same cell we have
 \begin{equation} \label{e28}
|z_j-z'_j|\le
\cases{{\e\over \log\log\nu}&,\ $1\le j\le\nu^{1/2}$,\cr
       \e                   &, $\nu^{1/2} <j\le
\nu$.}
\end{equation}
Thus, every cell is a product of two cubes of different
size and dimension. The necessary number of cells $M(\e)$ is bounded as follows,
$$
M(\e)\le \left({\log\log\nu\over \e}\right)^{[\nu^{1/2}]} \e^{-(\nu-[\nu^{1/2}])} = (1/\e)^\nu
(\log\log\nu)^{[\nu^{1/2}]}.
$$
Let us now evaluate the distance $\rho(z,z')$ for $z,z'$ satisfying
(\ref{e28}). By (\ref{e27}) we have
$$
\rho(z,z')\le C N^{1/2-\s} \D \left\{ \rho_1+\rho_2+\rho_3 \right\},
$$
where
$$
  \rho_1= \sum_{j=1}^\nu |z_{j} - z'_{j}| \sum_{b=2}^\infty b^{1+H} \ p_{j}^{-b},
$$
$$
  \rho_2 =\sum_{\nu^{1/2}<j\le \nu}  |z_{j} - z'_{j}| p_{j}^{-1},
$$
$$
   \rho_3  =\sum_{j\le \nu^{1/2}}  |z_{j} - z'_{j}| p_{j}^{-1}.
$$
For any $j\ge 1$ we have
\begin{equation} \label{e29}
  \sum_{b=2}^\infty b^{1+H} \ p_{j}^{-b} =
  \sum_{b=2}^\infty b^{1+H} \ ({2\over p_{j}})^{b} 2^{-b}
  \le  ({2\over p_{j}})^{2} \sum_{b=2}^\infty b^{1+H} \ 2^{-b}
  = C  p_{j}^{-2}.
\end{equation}
Hence,
$$
\rho_1\le \left(\sum_{j=1}^\nu C  p_{j}^{-2}\right)\ \ \max_{j\le \nu}  |z_{j} - z'_{j}| \le C \e.
$$
Similarly,
$$ \rho_2\le \left(\sum_{\nu^{1/2}<j\le\nu} p_{j}^{-1}\right)\ \ \max_{\nu^{1/2}< j\le \nu}  |z_{j} - z'_{j}|
\le C\
\left(\sum_{\nu^{1/2}<j\le\nu}  (j\log j)^{-1}\right)\ \e
$$
$$
\le C \int_{\nu^{1/2}}^{\nu} {du\over u\log u}\ \e
 = C (\log\log \nu-\log({\log\nu\over 2})) \ \e = C  (\log2) \, \e.
$$
Finally,
$$
\rho_3\le \left(\sum_{j=1}^\nu p_{j}^{-1}\right)\ \ \max_{j\le \nu^{1/2}} |z_{j} - z'_{j}|
\le C \left(\sum_{j=1}^\nu (j\log j)^{-1}\right)\ {\e\over \log\log\nu} \le C\ \e.
$$
By summing up three estimates, we have
$\rho(z,z')\le C N^{1/2-\s} \D \e$ which enables the evaluation of the metric entropy.

Let $\NNN\left(\T^\nu,\rho,u\right)$ be the minimal number of balls of
radius $u$ that cover the space $(\T^\nu,d)$. We have
$$
\log\NNN\left(\T^\nu,\rho, C N^{1/2-\s}\ \D \e\right) \le \log M(\e)
\le \nu |\log\e|+ \nu^{1/2}\cdot \log\log\log\nu  .
$$
Observe also that
\begin{equation} \label{e210}
\|\Upsilon (\uz)\|_2\le C N^{1/2-\s} \D ,\quad \uz\in
\T^{\nu}.
\end{equation}
Hence,
$D:=diam(\T^\nu,\rho)\le C N^{{1\over 2}-\s}\D  $,
and by the classical Dudley's  entropy theorem (see \cite{L}, Theorem 1 p.179),
for any fixed $\uz\in \T^\nu$
\begin{eqnarray*}
\E\,  \sup_{\uz'\in T^\nu} |\Upsilon (\uz')-\Upsilon(\uz)|
&\le& C\int_0^D [\log\NNN(\T^\nu,\rho,u)]^{1/2}du
\cr
&\le&  C \int_0^{C N^{1/2-\s}\D }
 [\log\NNN(\T^\nu,\rho,u)]^{1/2}du
\cr
&=& C N^{1/2-\s} \D
 \int_0^1 [\log\NNN(\T^\nu,\rho,C N^{1/2-\s}\e)]^{1/2}d\e
\cr
&\le& C N^{1/2-\s} \D \int_0^1\left[ \nu |\log\e|+ \log\log\log\nu \cdot \nu^{1/2} \right]^{1/2} d\e
\cr
&\le&    C N^{1/2-\s} \D \nu^{1/2}.
\end{eqnarray*}
Using again (\ref{e210}), we have
\begin{equation} \label{e211}
\E\,  \sup_{\uz'\in  T^\nu} |\Upsilon (\uz')|
\le    C
N^{1/2-\s} \D  \nu^{1/2}.
\end{equation}
The final stage of the proof provides the optimal choice of the parameter $\nu$
balancing the quantities (\ref{e24}), (\ref{e25}), and (\ref{e211}).
In the first version, we brutally replace by $\D$ the quantities $\td_2(N/p_\nu)$,
resp. $\td_1(N/p_\nu)$ in (\ref{e24}) and (\ref{e25}). By taking any $\nu$ in the fairly vast range
\[
\frac {\log^2 N}{\log^2\log N} \le \nu \le \frac {N}{\log^2 N}\ ,
\]
we conclude that the main term is (\ref{e24}) and thus obtain the
first result of the theorem.

Moreover, assuming (\ref{extra}) to be verified, we can choose
$\nu=N^h$, where $h=h(b)=\frac {1} {2(b+1)}$ satisfies $h<1$. Due to the
choice of $h$, we have
\begin{equation} \label{comp}
b+h/2-1/2 = b((1-h)-h/2= \frac{4b^2+2b-1}{4(b+1)} < 0,
\end{equation}
since $b<(\sqrt{5}-1)/4$.
The estimate (\ref{e24_ex}) already provides a good order. We now take
care of the two remaining bounds. For (\ref{e25}) we have
\begin{eqnarray*}
\frac{C N^{1-\s}\td_1(N/p_\nu)}{\nu^{1/2}\log\nu} &\le&
\frac{C N^{1-\s}}{\log N} \ \cdot \
\frac{\td_2(N/p_\nu)}{\nu^{1/2}}
\cr
&\le&
\frac{C N^{1-\s}}{\log N} \ \cdot \
\frac{\td_2(C N/\nu)}{\nu^{1/2}}
\le
\frac{C N^{1-\s}}{\log N} \ \cdot \
N^{b(1-h)-h/2},
\end{eqnarray*}
which is good by  (\ref{comp}).
Similarly, for  (\ref{e211}) we have
\[
C N^{1/2-\s}\td_2(N)\nu^{1/2} \le \frac{C N^{1-\s}}{\log N}
\left( N^{b+h/2-1/2} \cdot \log N
\right),
\]
which is also good by  (\ref{comp}).
Therefore, we obtain the second result of the theorem.

Finally, let us give the necessary modifications for the exceptional case $\sigma=1/2$.
The first deviation of the proof occurs at (\ref{abel}) where we have this time
\begin{equation} \label{abel_2}
 \sum_{m\le M} d(m)^2 m^{-1} \le C \td_2^2(M) \log M.
\end{equation}
After substitution in (\ref{e21}) this yields
\begin{equation}  \label{e21_2}
\sum_{ \nu <j\le \tau} \sum_{n\in E_j}  d(n)^2 n^{-1} (\a_j -\a'_j)^2
\le
C \sum_{\nu<j\le\tau} (\a_j -\a'_j)^2  {\td _2^2(N/p_j)  \log_+(N/p_j) \over p_j},
\end{equation}
Furthermore, at the place of (\ref{e23a}) we have
\begin{equation} \label{e23a_2}
K_m= {d(m)\over m^{1/2}} \ \left(\sum_{{\nu\le j\le\tau \atop mp_j\le N}} p_j^{-1}\right)^{1/2}
\le
C\ {d(m)\over m^{1/2}} \ [\log\log (N/m)]^{1/2}
\ .
\end{equation}
which results in
\begin{eqnarray}
\sum_{m\le N/p_\nu} K_m  &\le& C\  \sum_{m\le N/p_\nu}
 {d(m)\over m^{1/2}} \ [\log\log (N/m)]^{1/2}
\cr \label{e23b_2}
&\le&
{C \ N^{1/2} \td_1(N/p_\nu) (\log\log \nu)^{1/2} \over \nu^{1/2} (\log\nu)^{1/2}}\ .
\end{eqnarray}
at the place of (\ref{e23b}) and (\ref{e25}).
Next, at the place of (\ref{e24}) we have
\begin{eqnarray} \label{e24_2}
 \E\, \sup_{\g\in \G}|Y'(\g)|
 &\le& C \, \sum_{\nu\le j\le \tau} {\td _2(N/p_j)  \log_+^{1/2}(N/p_j) \over p_j^{1/2}}
 \cr
 &\le& C \, \td _2(N/p_\nu) \, \sum_{\nu\le j\le \tau} { \log_+^{1/2}(N/p_j) \over p_j^{1/2}}
 \cr
 &\le& {C\ N^{1/2}  \td_2(N/p_\nu)  \over \log N}
 \ ,
\end{eqnarray}
where at the last step we applied Lemma \ref{lem2} to the logarithmic function. Moreover,
under assumption (\ref{extra}) one applies Lemma \ref{lem2} to the function $f(r)=r^b\log_+^{1/2}r$,
and obtains a bound
\begin{equation} \label{e24_ex_2}
\E\, \sup_{\g\in \G}|Y'(\g)|
\le {C\ N^{1/2} \over \log N}\ ,
\end{equation}
which perfectly matches (\ref{e24_ex}).

The evaluation of the process $Q_1$ goes along the same lines as before until we have to use (\ref{abel_2})
thus arriving to
\begin{equation} \label{e27_2}
 \rho(z,z')\le C (\log N)^{1/2} \D \left\{
   \sum_{j=1}^\nu  |z_{j} - z'_{j}|
   \sum_{b=1}^\infty  b^{1+H} \  p_{j}^{-b} \right\}.
\end{equation}
at the place of  (\ref{e27}). One observes that only the constant is different :
$(\log N)^{1/2}$ stands for $N^{1/2-\s}$.
By  carrying  on this constant along the rest of the calculation, we arrive at
\begin{equation} \label{e211_2}
\E\,  \sup_{\uz'\in  T^\nu} |\Upsilon (\uz')|
\le    C  (\log N)^{1/2}\, \D  \nu^{1/2}.
\end{equation}
 at the place of  (\ref{e211}).
The concluding part of the proof which deals with optimal balance of three expressions
does not require any substantial change since the most sensible expression
(\ref{e24}), resp. (\ref{e24_2}) is the same in both cases.

 \cqfd


\bigskip
\section{Other results}
\medskip\par Fix some positive
integer $\tau\le \pi(N)$, and recall that $p_1<p_2<\ldots$
is the sequence of primes. Put
$$
{\cal E}_\tau= {\cal E}_\tau(N)=
\big\{ 2\le
n\le N : P^+(n)\le p_\tau\big\}.
$$
Note that for $\m=\pi(N)$ we have
${\cal E}_\m=\{2,\ldots,N\}$.
The ${\cal E}_\tau$-based Dirichlet polynomials were
considered in \cite{Q3} and \cite{LW}. In this section, we will
establish the theorem below extending the main result  of
\cite{LW} (Theorem 1.1) to the boundary case $\s=1/2$.

\begin{thm}
a) Upper bound. Then there exists a
constant $C $ such that for any integer $N\ge 2$
it is true that
 $$
\E \sup_{t \in \R} \big| \sum_{n\in {\cal E}_\tau}\e_n  n^{ -1/2 -it} \big|
\le
\cases{
 C  \left( {\tau\over \log\tau} \log_+({N\over p_\tau}) \right)^{1/2}
&,\
if \ $(N\log\log N)^{1/2} \le \tau\le {N\over \log N},$  \cr
& \cr
 C  \left(N \ \log\log N \right)^{1/4}
   &,\
if \ $ {(N\log\log N)^{1/2}\over \log N}\le  \tau \le (N\log\log N)^{1/2},$
\cr
& \cr
 C   \left(\log N \ \tau \right)^{1/2}
   &,\
if \  $1\le \tau \le {(N\log\log N)^{1/2}\over \log N} \ .$
} $$

b) Lower bound. There exists a constant
$c$ such that for every $N\ge 2$,
$$
\E\, \sup_{t \in \R} \big| \sum_{n\in {\cal E}_\tau}\e_n  n^{ -1/2 -it}\big|
\ge c  \ \left( {\tau\over \log\tau} \
\min\left\{
\log_+\left({N\over p_\tau}\right); \log p_{\tau/2}
\right\}\right)^{1/2}  .
$$
\end{thm}\bigskip

{\bf Proof of the upper bound}.
We start exactly as in  the proof of Theorem \ref{t1}  until (\ref{e21}) where one has
\begin{equation} \label{e32}
   \sum_{ \nu <j\le \tau} \sum_{n\in E_j} n^{-2\s}  (\a_j -\a'_j)^2
   \le
   \sum_{\nu<j\le\tau} (\a_j -\a'_j)^2p_j^{-2\s}  \sum_{m\le N/p_j} m^{-2\s}.
\end{equation}
When $\s=1/2$, it follows that
\begin{equation} \label{e33}
    \sum_{ \nu <j\le \tau} \sum_{n\in E_j} n^{-2\s}  (\a_j -\a'_j)^2
    \le
    C_\s \sum_{\nu<j\le\tau} (\a_j -\a'_j)^2 p_j^{-1}\log_+\big(N/p_j\big).
\end{equation}
Next, in (\ref{e22}) one obtains
\begin{eqnarray}
\sum_{\nu <j\le \tau} \sum_{n\in E_j}
{ (\b_{{n\over p_j}}- \b'_{{n\over p_j}})^2\over n^{2\s}}
&\le &
\sum_{m \le N/p_\nu} (\b_m-\b'_m)^2
 \big(\sum_{{\nu<j\le \tau\atop  mp_j\le N}}{1\over (mp_j)^{2\s}}
 \big)
\cr  \label{e34}
&:=&
\sum_{m \le N/p_\nu} K_m^2 (\b_m-\b'_m)^2.
 \end{eqnarray}
while  with $\s=1/2$ we have
$$K_m^2 = \sum_{{\nu<j\le \tau \atop  mp_j\le N}} (mp_j)^{-1}
= m^{-1} \sum_{{\nu<j\le \tau \atop  mp_j\le N}} p_j^{-1}.
$$
The upper summation border is different in two cases ($p_j\le p_\tau${\it vs} $p_j\le N/m$).
Therefore, we distinguish two cases.

1) $m\le N/p_\tau$. Here, $N/m\ge p_\tau$ and the border $j\le \tau$ is crucial. We obtain
\begin{eqnarray*}
    K_m^2 &=& m^{-1} \sum_{\nu<j\le \tau } p_j^{-1}
    \le  m^{-1} \sum_{1\le j\le \tau } p_j^{-1}
    \cr
    &\le&
    C \ m^{-1} \sum_{1\le j\le \tau }(j \log j)^{-1} \le   C\ m^{-1} \log\log\tau .
\end{eqnarray*}
It follows that
$$      \sum_{m\le N/p_\tau} K_m
   \le  C \sum_{m\le N/p_\tau} m^{-1/2} (\log\log\tau)^{1/2}
   \le    C (N/p_\tau)^{1/2}(\log\log\tau)^{1/2}.
$$

2) $N/p_\tau \le m \le N/p_\nu$. Here, the border $p_j\le N/m$ is crucial.
We choose a positive integer $k$ such that $N/p_k\sim m$ and obtain
$$
K_m^2 \le m^{-1} \sum_{j\le k} p_j^{-1}
      \le C m^{-1} \sum_{1\le j\le k } (j \log j)^{-1}
      \le C\ m^{-1} \log\log k
$$
$$  \sim  C\ m^{-1} \log\log p_k \sim  C\ m^{-1} \log\log (N/m).
$$
It follows
$$
K_m \le  C\ m^{-1/2} [\log\log (N/m)]^{1/2}.
$$
Hence,
\begin{eqnarray*}
\sum_{m\le N/p_\nu} K_m &\le &
C \sum_{m\le N/p_\nu} m^{-1/2} [\log\log (N/m)]^{1/2}\cr
& =&
C N^{-1/2} \sum_{m\le N/p_\nu} (m/N)^{-1/2} [\log\log (N/m)]^{1/2}
\cr &\le &
C N^{1/2} \int_{0}^{1/p_\nu} u^{-1/2} [\log\log (1/u)]^{1/2} du
\cr &\le&
C \left( N\over p_\nu \right)^{1/2} [\log\log p_\nu]^{1/2}.
\end{eqnarray*}
As in \cite{LW} ,   we define a second Gaussian process by
putting for all $\g\in \G$
\begin{eqnarray*}
Y(\g)& =&
\sum_{\nu<j\le\tau}
\big( p_j^{-1}\log(N/p_j)  \big)^{1/2}   \a_j\xi'_j
+
\sum_{m\le N/p_\nu} K_m \ \b_m\xi''_m
 :=
 \ Y'_\g + Y''_\g ,
\end{eqnarray*}
where $\xi'_i  $, $\xi''_j$  are independent
${\cal N}(0,1)$ random
variables. It follows from (\ref{e33})
  and
(\ref{e34}) that for some suitable constant
$C_\s$, one has the
comparison relations: for all $\g, \g'\in \G$,
$$\|X_\g-X_{\g'}\|_2\le C_\s \|Y_\g-Y_{\g'}\|_2.
$$
Next, by virtue of the Slepian comparison lemma,
 since
$X_0=Y_0=0$,
we have
$$\E\, \sup_{\g\in \G}     |X_\g|\le
2 \E\,  \sup_{\g\in \G}      X_\g \le
2 C_\s \E\, \sup_{\g\in \G} Y_\g \le
2 C_\s \E\, \sup_{\g\in \G} |Y_\g|.
$$
It remains to evaluate the supremum  of $Y$. First of all,
 $$\E\, \sup_{\g\in \G}|Y'(\g)|
\le  \sum_{\nu<j\le\tau}   p_j^{-1/2}[\log(N/p_j)]^{1/2} .
$$
We have
\begin{eqnarray*}
   \sum_{\nu<j\le\tau}   p_j^{-1/2} [\log(N/p_j)]^{1/2}  &\le&
   C \sum_{1<j\le\tau}   (j \log j)^{-1/2} [\log(N/p_j)]^{1/2}
   \cr &\le &
   C \left({\tau\over \log\tau}\right)^{1/2}\, [\log(N/p_\tau)]^{1/2}.
\end{eqnarray*}
thus
\begin{equation} \label{e35}
   \E\, \sup_{\g\in \G}|Y'(\g)| \le
   C\ {\tau^{1/2} \over (\log\tau)^{1/2}} \ [\log(N/p_\tau)]^{1/2} .
\end{equation}
To control the supremum of $Y''$, we use our estimates for the sums of $K_m$ and
write that
\begin{eqnarray*}
  \E\, \sup_{\g\in \G} |Y''(\g)| &\le& \sum_{m\le  N/p_\nu} K_m
  \cr &\le&
   C \left({N \over p_\tau}\right)^{1/2}(\log\log\tau)^{1/2}
 +C \left({N\over p_\nu}\right)^{1/2} [\log\log p_\nu]^{1/2}.
\end{eqnarray*}
Since
$$ {\log\log\tau \over p_\tau} \le  {\log\log \nu \over p_\nu} \sim
{\log\log p_\nu \over p_\nu},
$$
only the second term is important and by
$$\log\log p_\nu\sim\log\log \nu$$
we obtain
\begin{equation} \label{e36}
 \E\, \sup_{\g\in \G} |Y''(\g)| \le
 C \left({N\over p_\nu}  \ \log\log\nu   \right)^{1/2}.
\end{equation}
Next, as in \cite{LW} , we  turn to the supremum of $Q_1 $.
Towards this aim, introduce the auxiliary Gaussian process
$$  \Upsilon
(\uz) =
\sum_{  P^+(n) \le p_\nu }   n^{-\s}
\big\{\t_n \cos 2\pi \langle
\ua(n),\uz\rangle +
\t_n'\sin 2\pi \langle \ua(n),\uz\rangle \big\}
,\qq
\uz\in \T^{\nu },
$$
where $\t_i$, $\t'_j$  are independent ${\cal N}(0,1)$ random variables. By symmetrization,
 $\displaystyle{\E\, \sup_{\uz \in \T^{\nu }}\big|Q_1(\uz)\big|\le
\sqrt{8\pi} \E\, \sup_{\uz \in  \T^{\nu }}\big|\Upsilon (\uz)\big|}$, so
that we are again led to evaluating the supremum of a real valued
Gaussian process. For  $\uz, \uz' \in \T^{\nu}$ put
\[
d (\uz,\uz') :=\, \big\|\Upsilon(\uz)-\Upsilon(\uz)\big\|_2,
\]
and observe that
\begin{eqnarray}
\,d(\uz,\uz')^2  &= &4 \sum_{n: P^+(n)\le p_\nu}
  {1\over n^{2\s}}
  \sin^2(\pi \langle\ua(n),\uz -\uz'\rangle)
\cr &\le &
  4\pi^2\
\sum_{n: P^+(n) \le p_\nu}  {1\over n^{2\s}}
    |\langle\ua(n),\uz -\uz'\rangle|^2
\cr
  & \le &
4\pi^2\ \sum_{n: P^+(n)\le p_\nu}   n^{-2\s}
\left[\sum_{j=1}^\nu a_j(n) |z_j - z'_j|\right]^2
\cr
&  = &
4\pi^2\ \sum_{n: P^+(n)\le p_\nu}   \sum_{j_1,j_2=1}^\nu a_{j_1}(n)
a_{j_2}(n) |z_{j_1} - z'_{j_1}|\ |z_{j_2} - z'_{j_2}| n^{-2\s}
\cr
& = &
4\pi^2\ \sum_{j_1,j_2=1}^\nu \sum_{n: P^+(n)\le p_\nu}
   a_{j_1}(n)
a_{j_2}(n)  |z_{j_1} - z'_{j_1}|\ |z_{j_2} - z'_{j_2}| n^{-2\s}
\cr
& \le &
4\pi^2\ \sum_{j_1,j_2=1}^\nu  |z_{j_1} - z'_{j_1}|\ |z_{j_2} -z'_{j_2}|
\sum_{b_1,b_2=1}^\infty b_1 b_2
\sum_{n\le N, a_{j_1}(n)=b_1,
a_{j_2}(n)=b_2}  n^{-2\s}
\cr
&\le &  4\pi^2  \sum_{j_1,j_2=1}^\nu |z_{j_1} - z'_{j_1}|\ |z_{j_2} - z'_{j_2}|
\sum_{b_1,b_2=1}^\infty  b_1 b_2 p_{j_1}^{-2 b_1\s} p_{j_2}^{-2 b_2\s}
\sum_{{k\le N p_{j_1}^{-b_1} p_{j_2}^{-b_2} \atop P^+(k)\le p_\nu}} k^{-2\s}.
\cr
&{}&  \label{e37}
\end{eqnarray}

In the case $\s=1/2$ we obtain
$$
d(\uz,\uz')^2 \le
4\pi^2\ \sum_{j_1,j_2=1}^\nu |z_{j_1} - z'_{j_1}|\ |z_{j_2} - z'_{j_2}|
\sum_{b_1,b_2=1}^\infty  b_1 b_2 p_{j_1}^{-b_1} p_{j_2}^{-b_2}
\log_+\left({N \over p_{j_1}^{b_1} p_{j_2}^{b_2} } \right).
$$
We ignore the denominator in the logarithm and move $\log N$ ahead:
\begin{equation} \label{e38}
d(\uz,\uz')^2 \le C\ \log N\
\left\{ \sum_{j=1}^\nu  |z_{j} - z'_{j}| \sum_{b=1}^\infty b \ p_{j}^{-b} \right\}^2,
\end{equation}
or, equivalently
\begin{equation} \label{e39}
d(\uz,\uz')^ \le C\ (\log N)^{1/2}\
 \sum_{j=1}^\nu  |z_{j} - z'_{j}| \sum_{b=1}^\infty b \ p_{j}^{-b}.
\end{equation}

The subsequent estimates are identical to those of \cite{LW}  with
$N^{1-2\s}$ being replaced by $\log N$
and $N^{1/2-\s}$  by $(\log N)^{1/2}$ until we arrive at
\begin{equation} \label{e40}
    \E\,  \sup_{\uz'\in  T^\nu} |\Upsilon (\uz')|
    \le    C_\s   (\log N)^{1/2} \nu^{1/2}.
\end{equation}
The final stage of the proof provides the optimal choice of
the parameter $\nu$ balancing the quantities  (\ref{e35}),
(\ref{e36}), and (\ref{e40}). As suggests the Theorem's claim, we consider three cases.

{\bf Case 1.}\ $(N\log\log N)^{1/2}  \le  \tau \le {N\over \log N}.$

In this case we choose
\begin{equation} \label{e41}
    \nu= {(N \, \log\log N)^{1/2} \over\log N}
\end{equation}
thus balancing (\ref{e36}) and (\ref{e40}). Both terms yield then $(N\, \log\log N)^{1/4}$
which is dominated  in this zone by the constant (\ref{e35}).
>From (\ref{e35}) we obtain the bound
$    C_{\s} \left( {\tau\over \log\tau}\ \log_+(N/p_\tau) \right)^{1/2}.
$
The correctness condition $\nu\le \tau$ is obvious.
\medskip

{\bf Case 2.}\ ${(N\log\log N)^{1/2} \over\log N} \le  \tau \le (N\log\log N)^{1/2}.$

In this case we still choose $\nu$ from (\ref{e41}) thus balancing (\ref{e36}) and (\ref{e40}) and
getting the bound $(N\, \log\log N)^{1/4}$. The difference is that in this range
the constant term (\ref{e35}) is negligible. It follows that our total bound is $(N\, \log\log N)^{1/4}$.
The correctness condition $\nu\le \tau$ is still obvious for
the range under consideration.
\medskip

{\bf Case 3.}\   $1\le \tau \le {(N\log\log N)^{1/2} \over\log N}.$
Here we just set $\nu=\tau$. It means that we do not
need the splitting of the polynomial in two parts. Formally,
the quantities (\ref{e35}) and (\ref{e36}) are not necessary and we
obtain the bound $C_{\s} \left(\log N \ \tau \right)^{1/2}$
directly from (\ref{e40}).
\smallskip\par\noi
The upper bound is proved completely.
\bigskip

{\bf Proof of the lower bound}.  We shall first collect some   auxiliary  useful results, valid for 
general Dirichlet polynomials. Next we will apply them to the considered  setting.
 Let $\ud =\{d_n, n\ge 1\}$ be a  sequence of reals. Recall
that by (\ref{e11}) we have
$$
\sup_{t\in \R}\big|\sum_{j=1}^\tau \sum_{n\in E_j} d_n\e_n n^{ -\s - it}\big|
=\sup_{\uz \in \T^\tau}\big|Q(\uz)\big|.
$$
where
$$ Q(\uz)= \sum_{j=1}^\tau \sum_{n\in E_j}
   d_n \e_n n^{-\s} e^{2i\pi\langle \ua(n),\uz\rangle}.
$$
\noi
Consider the subset $\zz$ of $\T^\tau$ defined by
$$\zz=\Big\{ \uz=\{z_j,
1\le j\le \tau\} : \ \hbox{$z_j=0$,
if $j\le \tau/2$,\
and
\ $z_j\in\{0,1/2\}$, if $j\in(\tau/2,\tau]$} \Big\} .
$$

Observe that the imaginary part of $Q$ vanishes on $\zz$, since for any
$\uz\in \zz$ and any $n$ it is true that
$$ e^{2i\pi\langle \ua(n),\uz\rangle} =
\cos(2\pi\langle
\ua(n),\uz\rangle) = (-1)^{2\langle
\ua(n),\uz\rangle}.
$$
Hence, $Q$ takes the following simple form on $\zz$
$$
 Q(\uz) = \sum_{\tau/2<j\le \tau}  \sum_{n\in E_j} d_n \e_n n^{-\s}
{(-1)}^{2\langle \ua(n),\uz\rangle }.
$$
This is no longer a trigonometric polynomial, but simply a finite
rank Rademacher process.

For $j\in(\tau/2,\tau]$ define
$$ {\cal L}_j=\Big\{n=p_j \, \tilde n\ : \
  \tilde n\le {N\over p_{j}}\ \hbox{and}\
P^+(\tilde n)\le p_{\tau/2}\Big\}.
$$

Since
$$
 E_j \supset{\cal L}_j, \q j=1,\ldots \tau,
$$
the sets ${\cal L}_j$ are pairwise disjoint.

Put for $z\in \zz$,
$$
 Q'(\uz)= \sum_{  \tau/2<j\le \tau } \sum_{n\in  {\cal L}_j}
 d_n \e_n n^{-\s} {(-1)}^{2\langle \ua(n),\uz\rangle }.
 $$
We now recall a useful fact. (\cite{LW} Lemma 3.1)

\begin{lem}   Let  $X=\{X_z, z\in Z\}$ and
$Y=\{Y_z,
z\in Z\}$ be two finite sets of random variables defined on
a common probability space. We assume that $X$ and $Y$ are
independent and that the random variables $Y_z$ are all centered.
Then
$$ \E\sup_{z\in
Z}|X_z + Y_z| \ge  \E\sup_{z\in Z}|X_z  |.$$
\end{lem}

Clearly, since
$\{Q(\uz)-Q'(\uz),\uz\in \zz\}$ and $\{Q'(\uz),\uz\in \zz\}$ are independent,
$$
\E \sup_{\uz\in \zz} |Q(\uz)| \ge
\E \sup_{\uz\in \zz} \left|Q'(\uz) \right|.
$$
We now proceed to a direct evaluation of $Q'(\uz)$ by proving
\medskip

\begin{prop}   There exists a universal constant $c$ such that for any system of coefficients
$(d_n)$
$$c\ \sum_{  \tau/2<j\le \tau }
\big|\sum_{n\in  {\cal L}_j}
d_n^2
\big|^{1/2} \le \E\, \sup_{\uz\in \zz}
\left|Q'(\uz)
\right|
\le
\sum_{  \tau/2<j\le \tau }
\big|\sum_{n\in  {\cal L}_j} d_n^2 \big|^{1/2}.
$$
 \end{prop}
\begin{cor}  If $(d_n$) is a multiplicative system (namely $d_{nm}=d_nd_m$ if $n,m$ are coprimes), we have
$$ \E\, \sup_{\uz\in \zz} \left|Q'(\uz) \right|\ge
c\ \sum_{  \tau/2<j\le \tau } d_{p_j}\Big( \sum_{\tilde n\le N/p_j\atop
  P^+(\tilde n) \le p_{\tau /2}}  d^2_{\tilde n}\Big)^{1/2}.
$$
\end{cor}

\noi When $\s=1/2$, we apply Corollary 3.3 with $d_n=n^{-1/2}$ which yields
\begin{eqnarray*}
\E\, \sup_{t \in \R} \big| \sum_{n\in {\cal E}_\tau}\e_n  n^{ -1/2 -it}\big|
&\ge&
c\ \sum_{  \tau/2<j\le \tau } p_j^{-1/2}
    \Big( \sum_{n\le N/p_j\atop P^+(n) \le p_{\tau /2}}
    n^{-1} \Big)^{1/2}
\cr & \ge &c\ {\tau\over 2}\ p_\tau^{-1/2} \
\Big( \sum_{n\le \min\{ N/p_\tau; p_{\tau /2}\}}
    n^{-1} \Big)^{1/2}
\cr  &\ge  &\left({c\ \tau\over \log\tau}\right)^{1/2}  \
\Big( \log_+ \min\{ N/p_\tau; p_{\tau /2}\}  \Big)^{1/2}.
\end{eqnarray*}
\bigskip

{\bf Acknowledgements.}\   The work of the first mentioned author was supported by
grants RFBR 05-01-00911 and NSh-422.2006.1 "Leading scientific schools of
Russia".  He is also grateful for hospitality to the L.Pasteur University (Strasburg)
where this research has been done.
\medskip


{
}

\noi   {\phh Mikhail   Lifshits,  St.Petersburg State University,
Department of   Mathematics and Mechanics, 198504, Bibliotechnaya pl, 2,
Stary Peterhof, Russia.
\par\noindent
E-mail:\ \tt lifts@mail.rcom.ru}

\noi {\phh Michel  Weber, \noi  Math\'ematique (IRMA),
Universit\'e Louis-Pasteur et C.N.R.S.,   7  rue Ren\'e Descartes,
67084 Strasbourg Cedex, France.
\par\noindent
E-mail: \  \tt weber@math.u-strasbg.fr}
\end{document}